\documentclass{article}
\usepackage{graphicx}
\usepackage{amsmath}
\usepackage{amsfonts}
\usepackage{amsthm}

\begin{document}

\title{Transformations  Integer Sequences 
\\And Pairing Functions}

\author{Boris Putievskiy}
\date{December 13, 2012}
\maketitle

\begin{abstract}
We propose several procedures for creating new families of integer sequences based on the method of Cantor diagonalization. Then we modify and generalize this method. The paper includes explicit formulas for most proposed families of integer sequences.

\end{abstract}

\tableofcontents

\section{Introduction}
There are well-known transformations of the integer sequences [1],[2],[3].
Let $\mathcal{A}$ be the set of integer sequences. We construct an array $\Omega:  \{ w(i,j) \} $, in which each column is an element of  $\mathcal{A}$.
We define the  sequence $\omega(n)$ $\in$ $\mathcal{A}$ as the list of the elements of the array $\Omega$
in the Cantor anti-diagonal order [4]. All sequences in the paper are numbered starting from 1. New sequence starts with:
\\$\omega(1)=w(1,1), 
\\\omega(2)=w(1,2), \quad \omega(3)=w(2,1), 
\\\omega(4)=w(1,3), \quad \omega(5)=w(2,2), \quad \omega(6)=w(3,1), $
\\. . .
\\The  sequence $\omega$ is a triangular array read by rows: row number $k$ contains $k$ elements. 
From [4],[5] the formula for calculating $n$ is:
\begin{equation}
\\n=\dfrac{(i+j-1)(i+j-2)}{2} + i
\end{equation} 
The property of (1) is that rearranging  $i$, $j$ we get the list of elements of the array $\Omega$ in the Cantor diagonal order.
\\The formula (1) is obtained by simple transformation of the classical formula:
\begin{equation}
\\z=\dfrac{(x+y)^{2}+3x+y}{2}
\end{equation} 
The characteristic of (2) is that the enumeration starts from pair $x=0$, $y=0$  to $z=0$.
\\
\\Using C. Kimberling formula from A002260, we get inverse functions for (1):
\\$i=n-\dfrac {t(t+1)}{2}, \quad j=\dfrac {t^2+3t+4}{2}-n$,
\\where $t=\lfloor \dfrac{\sqrt{8n-7}-1}{2} \rfloor $.
\\\\Throughout the paper, sequences will be referred to by their Annnnnn  numbers, as found in the On-Line Encyclopedia of Integer Sequences [6].
\\
\\We denote the set of integers numbers $\mathbb{Z}$, the set of the integer nonnegative numbers $\mathbb{Z}^{*}$, and the set of the integer positive numbers $\mathbb{Z}^{+}$. 
\section{Transformation of a Single Sequence}
\subsection{Replication}
{\bf Example 2.1.1} \\
Let $\alpha$ $\in$ $\mathcal{A}$: $a_{1}, a_{2}, a_{3},...$ We form the array $\Omega$ replicating the column $\alpha$. It follows that $w(i,j)=a_{i}$ for any $j$:
\\
\begin{equation} 
\begin{matrix}
a_1 & a_1& a_1 \quad.\quad.\quad. \\
a_2 & a_2& a_2 \quad.\quad.\quad. \\
a_3 & a_3& a_3 \quad.\quad.\quad. \\
\vdots& \vdots& \ddots \\
\end{matrix}
\end{equation} 
\\
\\and $\omega(n)$ is formed as:
\\$a_{1},
\\a_{1},a_{2},
\\a_{1},a_{2},a_{3},
\\ \quad .\quad .\quad.$
\\
The resulting $\omega$ is a fractal sequence [7]. 
\\
\\{\bf Definition 1.} A sequence $\beta$ is called a reluctant sequence of sequence $\alpha$, if $\beta$ is triangle array read by rows: row number $k$ coincides with first $k$ elements of the sequence $\alpha$.
\\\\Some of reluctant sequences are presented in OEIS A002260, A002262, A037126, A059268, A215026. 
\\
\\Formula for a reluctant sequence is:
\\$\omega(n)=a_{m}$, $m=n-\dfrac {t(t+1)}{2}$, where $t=\lfloor \dfrac{\sqrt{8n-7}-1}{2} \rfloor $.
\\
\\{\bf Example 2.1.2} \\
We apply the procedure of replication to the sequence $\omega$ from Example 2.1.1. 
\\Using double replication, we obtain $\omega^{'}$, wich is a reluctant sequence of the sequence  $\omega$ :
\\$a_{1},
\\a_{1},a_{1}, 
\\a_{1},a_{1},a_{2}, 
\\a_{1},a_{1},a_{2},a_{1},
\\a_{1},a_{1},a_{2},a_{1},a_{2},
\\a_{1},a_{1},a_{2},a_{1},a_{2},a_{3},
\\ \quad .\quad .\quad.$
\\$\omega^{'}(n)=a_{m}$,\quad $m=n^{'}-\dfrac {t^{'}(t^{'}+1)}{2}$, 
\\where 
\\$n^{'}=n-\dfrac {t(t+1)}{2}$, \quad $t^{'}=\lfloor \dfrac{\sqrt{8n^{'}-7}-1}{2} \rfloor $, \quad
$t=\lfloor \dfrac{\sqrt{8n-7}-1}{2} \rfloor $.
\\
\\{\bf Example 2.1.3} \\
We form the array $\Omega$ replicating the row $\alpha$, so $w(i,j)=a_{j}$ for any $i$. The array  $\omega(n)$ is the transposition of the array from Example 2.1.1:
\\$a_{1},
\\a_{2},a_{1},
\\a_{3},a_{2},a_{1},
\\\quad .\quad .\quad.$
\\
\\{\bf Definition 2.} A sequence $\beta$ is called a reverse reluctant sequence of sequence $\alpha$, if $\beta$ is a triangle array read by rows: row number $k$ lists first $k$ elements of the sequence $\alpha$ in reverse order.
\\\\There are some examples of reverse reluctant sequences - A004736, A025581, A104762, A104887, A130321.
\\
\\Formula for a reverse reluctant sequence is:
\\$\omega(n)=a_{m}$, $m=\dfrac {t^2+3t+4}{2}-n$, where $t=\lfloor \dfrac{\sqrt{8n-7}-1}{2} \rfloor $.
\subsection{Self-composition}
Let $\alpha$ $\in$ $\mathcal{A}$: $a_{1}, a_{2}, a_{3},...$ We will assume in this section $a(n) \in \mathbb{Z}^{+} $ for any $n$. The first column of the array $\Omega$ is the sequence $\alpha$: $a(n)$, the second column is $a_{a(n)}$, the third column is $a_{a_{a(n)}}$, and so on.
Using another notation, $w(i,j)= \underbrace{a(a(\ldots a(i)))}_{j} $.
\\We get new sequence  $\omega(n)$:
\\$a_{1}, \quad a_{a_{1}}, a_{2}, \quad a_{a_{a_{1}}}, a_{a_{2}}, a_{3}, \quad... $ 
\\
\\{\bf Example 2.2.1} \\
Let $m \in \mathbb{Z}^{+}$, $m > 1$ and $a(n)=m^n$, 
\\New sequence  $\omega(n)$ is:
\\ $m^{1}, \quad  m^{m^{1}}, m^{2}, \quad  m^{m^{m^{1}}},m^{m^{2}},m^{3}, \quad...$
\\
\\{\bf Example 2.2.2} \\
Suppose $a(n)$ is the sequence of prime numbers, then the first column of $\Omega$ is  the sequence A000040, the second colum of  $\Omega$ is the sequence A006450. The start of the sequence  $\omega(n)$ is:
\\2,  $\quad$   3, 3,    $\quad$  11, 5, 5,  $\quad$    5381, 31, 11, 7,    ...
\\
\\{\bf Example 2.2.3} \\
Let $\phi(n)$ be Euler's totient function. The first four columns of $\Omega$ are  sequences A000010, A010554, A049099, A049100. 
The start of the sequence  $\omega(n)$ is:
\\1,  $\quad$   1, 1,    $\quad$  1, 1, 2,  $\quad$    1, 1, 1, 2,  $\quad$   1, 1, 1, 1, 4,  $\quad$  1, 1, 1, 1, 2, 2, ...
\subsection{Function Acting on an Index}
Let sequence  $\alpha \in \mathcal{A}$ and function $f$: $\mathbb{Z}^{+}\times\mathbb{Z}^{+}$ $\rightarrow$ $\mathbb{Z}^{+}$. 
We create the array $\Omega$: $w(i,j)= a_{f(i,j)}$. 
\\
\\{\bf Example 2.3.1} \\
Let the first column of the the array $\Omega$ be the sequence $\alpha$: $a_{1}, a_{2}, a_{3},...$   
\\Every next column is formed from previous shifted by $k \in \mathbb{Z}^{*}$ elements. We get $f(i,j)=i+kj-k$.  
The  start of the sequence $\omega$ is:
\\$a_{1},\quad a_{k+1},a_{2},\quad a_{2k+1},a_{k+2},a_{3},\quad a_{3k+1},a_{2k+2},a_{k+3},a_{4},\quad ...$
\\
\\The sequence is the generalization of A002260 ($k$=0), A002024 ($k$=1), A128076 ($k$=2), A131914 ($k$=3). 
\\\\$\omega(n)=a_{m}$, 
where $m=k(t+1) + (k-1) \Bigr (\dfrac {t(t+1)}{2}-n \Bigl) $, 
\\\\  $t=\lfloor \dfrac{\sqrt{8n-7}-1}{2} \rfloor $. 
\\
\\Using rearranging  $i$ and $j$ in the formula $f(i,j)=ki+j-k$, we obtain the array $\Omega$, in which the first row is $\alpha$, 
every next row  is formed from previous shifted by $k$ elements.
\\
\\{\bf Example 2.3.2} \\
Let $\alpha$ $\in$ $\mathcal{A}$: $a_{1}, a_{2}, a_{3},...$,  $k \in \mathbb{Z}^{+}$, and $f(i,j)$= max $\{ki+j-k;i+kj-k\}$.
\\The  start of the sequence $\omega$ is:
\\$a_{1},\quad a_{k+1},a_{k+1},\quad a_{2k+1},a_{k+2},a_{2k+1},\quad  a_{3k+1},a_{2k+2},a_{2k+2}, a_{3k+1}...$
\\
\\The sequence is the generalization of A002024 ($k$=1), A204004 ($k$=2), A204008 ($k$=3). 
\\$\omega(n)=a_{m}$, 
\\where $m=k(t+1)+(k-1)$max $\{\dfrac{t(t+1)}{2}-n;n-\dfrac{t^2+3t+4}{2}\}$,
\\\\  $t=\lfloor \dfrac{\sqrt{8n-7}-1}{2} \rfloor $. 
\\
\\{\bf Example 2.3.3} \\
Let $\alpha$ $\in$ $\mathcal{A}$: $a_{1}, a_{2}, a_{3},...$, and $k \in \mathbb{Z}^{+}$. The first column of the array $\Omega$ is the sequence $\alpha$, in all columns with number $j \quad (j>1)$  the segment of the sequence $\alpha$ with the length of $j-1$: $<a_{k+j-2}, a_{k+j-3},...,a_{k}>$  shifts the sequence $a_{1},a_{2}, a_{3},...$. 
\\Therefore
\\$f(i,j) =\begin{cases}
i-j+1,&\text{if $i\geq j$},\\
j-i+k-1,&\text{if $ i<j $}.\\
\end{cases}$
\\
\\The  start of the sequence $\omega$ is:
\\$a_{1},\quad a_{k},a_{2},\quad a_{k+1},a_{1},a_{3},\quad  a_{k+2},a_{k},a_{2}, a_{4},\quad  a_{k+3},a_{k+1},a_{1},a_{3},a_{5},...$
\\
\\The sequence is the generalization of A143182 ($k$=2).
\\$\omega(n)=a_{m}$, 
\\where $m=|(t+1)^2 - 2n| + k\lfloor \dfrac{t^2+3t+2-2n}{t+1}\rfloor$,
\\\\  $t=\lfloor \dfrac{\sqrt{8n-7}-1}{2} \rfloor $. 
\section{Transformation of a Pair of Sequences}
We form the array $\Omega$ from two sequences  $\alpha, \beta \in \mathcal{A}$, $\alpha: a_{1}, a_{2}, a_{3},...$, 
\\$\beta: b_{1}, b_{2}, b_{3},...$ 
Let function $F$: $\mathbb{Z}\times\mathbb{Z}$ $\rightarrow$ $\mathbb{Z}$. We define the array $\Omega$ as 
\\$w(i,j)=F(a_i,b_j)$. The  start of the sequence $\omega$ is:
\\$F(a_1,b_1),\quad F(a_1,b_2), F(a_2,b_1), \quad F(a_1,b_3),F(a_2,b_2),F(a_3,b_1), \quad ... $
\\Using formulas for reluctant sequence and reverse reluctant sequence we get:
\\$\omega(n)=F(a_{m_1},b_{m_2})$, \quad $m_1=n-\dfrac {t(t+1)}{2}, \quad m_2=\dfrac {t^2+3t+4}{2}-n$,
\\where $t=\lfloor \dfrac{\sqrt{8n-7}-1}{2} \rfloor $.
\\We denote $\omega=C(\alpha,\beta,F)$.
\\
\\{\bf Example 3.1} \\
Let $w(i,j)=a_{i}*b_{j}$, $p$ and $q$ are primes, $\alpha$ is the sequence $p^{1},p^{2},p^{3},...$, 
\\$\beta$ is the sequence $q^{1},q^{2},q^{3},...$ 
so $\omega (n)$ are numbers of the form $p^{i}q^{j}$. The sequence is the development of A036561. 
\\
\\{\bf Example 3.2} \\
Let $k \in \mathbb{Z}^{+}$, $w(i,j)=a^{k}_{i}+b^{k}_{j}$,
hence $\omega (n)$ are numbers that can be represented as the sum of the $k$-th powers of elements of sequences $\alpha$ and $\beta$.
\\
\\{\bf Example 3.3} \\
We apply the transformation of a pair of sequences for enumeration n-tuple [8]. We denote concatenation of two integer positive numbers as $+\!\!\!+$. For example 12$+\!\!\!+$345 = 12345. Function $F$ is concatenation: $w(i,j)=a_{i}$+\!\!\!$+b_{j}$. 
\\Let $\eta ^{1}$ be 1,2,3,... We denote $\eta ^{d}(n)= C(\eta ^{d-1},\eta ^{1},$+\!\!\!+$)$.
\\The array $\Omega$ for $d=2$ is:
\begin{equation}
\begin{matrix}
11 & 12 & 13 \quad.\quad.\quad. \\
21 & 22 & 23 \quad.\quad.\quad. \\
31 & 32 & 33 \quad.\quad.\quad. \\
\vdots& \vdots& \ddots \\
\end{matrix}
\end{equation}
\\Starts of sequences $\eta ^{d}(n)$ are:
\\$d=2$ \quad 11,\quad 12,21, \quad 13,22,31, \quad 14,23,32,41,...\quad \quad (A066686)
\\$d=3$ \quad 111,\quad 112,121, \quad 113,122,211, \quad 114,123,212,131,...
\\$d=4$ \quad 1111,\quad 1112,1121, \quad 1113,1122,1211, \quad 1114,1123,1212,1131,...
\\$d=5$ \quad 11111,\quad 11112,11121, \quad 11113,11122,11211, \quad 11114,11123,11212,11131,...\\
\\
\\{\bf Example 3.4} \\
Composition of a pair of sequences. Let $a_{i},b_{i} \in \mathbb{Z^{+}}$ for all $i$. 
\\$F(a_i,b_j)= \underbrace{b(b(\ldots a_i))}_{j} $.
\\Start of the sequence $\omega(n)$ is: $b_{a_{1}}, \quad b_{b_{a_{1}}}, b_{a_{2}}, \quad b_{b_{b_{a_{1}}}},b_{b_{a_{2}}}, b_{a_{3}}, \quad ...$
\section{Transformation of Several Sequences}
We form the array $\Omega$ from $l$ ($l \in \mathbb{Z}^{+}, l>1$) sequences  $\alpha^{1}, \alpha^{2}, ... \alpha^{l}\in \mathcal{A}$. Elements of sequence $\alpha^{j}$ denote $a^{j}_{1}, a^{j}_{2}, a^{j}_{3}...$ Let functions $f,s$: $\mathbb{Z}^{+}\times\mathbb{Z}^{+}$ $\rightarrow$ $\mathbb{Z}^{+}$,   $1 \leq s\leq l$. 
We construct the array $\Omega$: $w(i,j)=a^{s(i,j)}_{f(i,j)}$.
\\
\\{\bf Example 4.1} \\
Let $\alpha^{1}$,$\alpha^{2}$, ... $\alpha^{l} \in \mathcal{A}$. The array $\Omega$ is the replication of $l$ columns $\alpha^{1}$,$\alpha^{2}$, ... $\alpha^{l}$:
\\$f(i,j)=i$,   $s(i,j) = 1+$Mod$(j-1,l)$.
\\The start of the sequence $\omega$ for $l=3$ is:
\\$a^{1}_{1},\quad a^{2}_{1},a^{1}_{2}, \quad a^{3}_{1},a^{2}_{2},a^{1}_{3}, \quad a^{1}_{1},a^{3}_{2},a^{2}_{3},a^{1}_{4},\quad a^{2}_{1},a^{1}_{2},a^{3}_{3},a^{2}_{4},a^{1}_{5},\quad ...$
\\Using the formula for the reluctant sequence, we get the number of the element $m$ of the sequence $\omega(n)=a^{r}_{m}$: 
\\$m=n-\dfrac {t(t+1)}{2}$, where $t=\lfloor \dfrac{\sqrt{8n-7}-1}{2} \rfloor $.  
\\\\Taking into account the formula from Example 2.1.3, we obtain the number of sequence $r$: 
\\$r=1+$Mod$(\dfrac{t^2+3t+4}{2}-n-1,l)$, 
\\where $t=\lfloor \dfrac{\sqrt{8n-7}-1}{2} \rfloor $.  
\\
\\{\bf Example 4.2} \\
Let $\alpha^{1}$,$\alpha^{2}$, ... $\alpha^{l} \in \mathcal{A}$. The array $\Omega$ is replication of ``a braid'' from $l$ sequences: 
\\$f(i,j)=i$ and $s(i,j)=1+$Mod$(i+j-2,l)$. 
\\\\The start of the sequence $\omega$ for $l=3$ is: 
\\$a^{1}_{1},\quad a^{2}_{1},a^{2}_{2}, \quad a^{3}_{1},a^{3}_{2},a^{3}_{3}, \quad a^{1}_{1},a^{1}_{2},a^{1}_{3},a^{1}_{4},\quad a^{2}_{1},a^{2}_{2},a^{2}_{3},a^{2}_{4},a^{2}_{5},\quad ...$
\\
\\Using the formula for the reluctant sequence, we get the number of the element $m$ of the sequence $\omega(n)=a^{r}_{m}$:
\\$m=n-\dfrac {t(t+1)}{2}$, where $t=\lfloor \dfrac{\sqrt{8n-7}-1}{2} \rfloor $.  
\\\\Taking into account the formula from Example 2.3.2 at $k=1$, we get the number of sequence $r$: 
\\$r=1+$Mod$(t,l)$, where $t=\lfloor \dfrac{\sqrt{8n-7}-1}{2} \rfloor $. 
\\
\\
{\bf Example 4.3} \\
Let $\alpha^{1}$,$\alpha^{2}$, ... $\alpha^{l} \in \mathcal{A}$ and $j \in \mathbb{Z}^{+},\quad j>1$. The segment of the sequence $\alpha^{l}$:
\\  $<a_{j-1}^{l}, a_{j-2}^{l},...,a_{1}^{l}>$  shifts in all columns with number $j$  the  replications of sequences $\alpha^{1}$,$\alpha^{2}$, ... $\alpha^{l-1}$. 
\\Thereby 
\\$f(i,j) =\begin{cases}
i-j+1,&\text{if $i\geq j$},\\
j-i,&\text{if $ i<j $}.\\
\end{cases}$
\\
$f(i,j)$ is the special case Example 2.3.3 at $k=1$. 
\\
\\$s(i,j) =\begin{cases}
1+$Mod$(j-1,l-1),&\text{if $ i \geq j $},\\
l,&\text{if $i < j$}.\\
\end{cases}$
\\\\The start of the sequence $\omega$ for $l=3$ is
\\$a^{1}_{1},\quad a^{3}_{1},a^{1}_{2}, \quad a^{3}_{2},a^{2}_{1},a^{1}_{3}, \quad a^{3}_{3},a^{3}_{1},a^{2}_{2},a^{1}_{4},\quad a^{3}_{4},a^{3}_{2},a^{1}_{1},a^{2}_{3},a^{1}_{5},\quad a^{3}_{5},a^{3}_{3},a^{3}_{1},a^{1}_{2},a^{2}_{4},a^{1}_{6},...$
\\\\Taking into account the formula from Example 2.3.3 at $k=1$, we obtain the number of the element $m$ of the sequence $\omega(n)=a^{r}_{m}$: 
\\$m=|(t+1)^2 - 2n| + \lfloor \dfrac{t^2+3t+2-2n}{t+1}\rfloor$, where $t=\lfloor \dfrac{\sqrt{8n-7}-1}{2} \rfloor $. 
\\Formula for number of sequence $r$ is:
\\$r=l+v\Bigr($Mod$(\dfrac{t^{2}+3t+4}{2}-n-1,l-1)-l+1\Bigl)$,
\\\\ where $t=\lfloor \dfrac{\sqrt{8n-7}-1}{2} \rfloor $, $v=\lfloor (2n-t(t+1) +1)/(t+3)  \rfloor $. 
\section{Other Pairing Functions}
So far we were using the classical Cantor anti-diagonal order. Next, we define several pairing functions $\mathbb{Z}^{+}\times\mathbb{Z}^{+}$ $\rightarrow$ $\mathbb{Z}^{+}$.
\subsection{Diagonal Enumeration}
We denote $<w_{1},_{j};w_{2},_{(j-1)};\quad...\quad w_{j},_{1}>$ a diagonal of the array $\Omega$.
Any permutation elements of the diagonal is the enumeration of the array $\Omega$.
\\
\\{\bf Example 5.1.1} \\
We form the permutation so that  neighboring diagonals traverse in opposite directions. Such sequences were presented previously without explicit formulas [9], A056023, A056011. 
\\
\\Obviously $i+j$ is constant for all elements of the diagonal. We do not change the order of the list for an odd $i+j$. We reverse the diagonal for an even $i+j$: 
\\\\$n=\dfrac{1}{2}\Bigr((i+j-1)(i+j-2) -((-1)^{i+j}-1)i + ((-1)^{i+j}+1)j \Bigl)$.
\\\\The result  of the enumeration of the array (4) is: 
\\${11},\quad {12},{21}, \quad {31},{22},{13}, \quad {14},{23},{32},{41},\quad {51},{42},{33},{24},{15},\quad {16},{25},{34},{43},{52},{61},...$
\\
\\{\bf Example 5.1.2} \\
We denote by $k=\lfloor \dfrac {j+1}{2} \rfloor$.
We form the permutation as a symmetrical movement from center to edges:
\\$<w_{k},_{k};w_{(k-1)},_{(k+1)};w_{(k+1)},_{(k-1)}; \quad...\quad w_{1},_{j}; w_{j},_{1}>$ for odd $j$,
\\$<w_{(k-1)},_{(k+1)};w_{(k+1)},_{(k-1)}; \quad...\quad w_{1},_{j}; w_{j},_{1}>$  for even $j$.
\\
\\\\$n =\begin{cases}
\dfrac{i(i+1) + (j-1)(2i+j-4)}{2},&\text{if $i \geq j$},\\\\
\dfrac{i(i+1) + (j-1)(2i+j-4)}{2} + 2(j-i)-1,&\text{if $i < j$}.
\end{cases}$
\\\\The result of the enumeration of the array (4) is:
\\${11},\quad {12},{21}, \quad {22},{13},{31}, \quad {23},{32},{14},{41},\quad {33},{24},{42},{15},{51},\quad {34},{43},{25},{52},{16},{61},...$
\\
\\
\\{\bf Example 5.1.3} \\
We denote by $k=\lfloor \dfrac {j+1}{2} \rfloor$. We form the permutation as a symmetrical movement from edges to center: 
\\$<w_{1},_{j}; w_{j},_{1}; w_{2},_{(j-1)};w_{(j-1)},_{2}; \quad...\quad w_{k},_{k}>$ for odd $j$,
\\$<w_{1},_{j}; w_{j},_{1}; w_{2},_{(j-1)};w_{(j-1)},_{2}; \quad...\quad w_{(k-1)},_{(k+1)};w_{(k+1)},_{(k-1)}>$  for even $j$.
Such sequences were presented previously without explicit formulas A064578, A194982.
\\\\$n =\begin{cases}
i(2i-1) + \dfrac{1}{2}(j-i)(3i+j-3),&\text{if $i \leq j$},\\\\
j(2j-1) + \dfrac{1}{2}(i-j)(3j+i-3) +1,&\text{if $i > j$}.
\end{cases}$
\\\\The result of the enumeration of the array (4) is:
\\${11},\quad {12},{21}, \quad {13},{31},{22}, \quad {14},{41},{23},{32}, \quad {15},{51},{24},{42},{33},\quad {16},{61},{25},{52},{34},{43},...$
\\
\\
\\{\bf Example 5.1.4} \\
We form the permutation as an asymmetrical movement from edges to center and inversely. Neighboring diagonals traverse in opposite directions: 
\\$<w_{j},_{1}; w_{2},_{(j-1)}; w_{(j-2)},_{3}; \quad...\quad w_{(j-1)},_{2}; w_{1},_{j}>$ for odd $j$,
\\$<w_{1},_{j}; w_{(j-1)},_{2}; w_{3},_{(j-2)}; \quad...\quad w_{2},_{(j-1)};w_{j},_{1}>$  for even $j$.
\\\\$n=\dfrac{1}{2}\Bigr((i+j-1)(i+j-2) +((-1)^{ \mathrm{max}\{i;j\}}+1)i - ((-1)^{\mathrm{max}\{i;j\}}-1)j \Bigl)$.
\\\\The result of the enumeration of the array (4) is:
\\${11},\quad {12},{21}, \quad {31},{22},{13}, \quad {14},{32},{23},{41}, \quad {51},{24},{33},{42},{15},\quad {16},{52},{34},{43},{25},{61},...$
\subsection{Angle Enumeration}
We consider squares of the array  $\Omega$ with vertices $w(1,1), w(1,i), w(i,i), w(i,1)$.
We list the array $\Omega$ according to the sides of squares from $w(1,i)$ to $w(i,i)$, then from $w(i,i)$ to $w(i,1)$.
The numeration was presented previously as A060734, A060736:
\\$n =\begin{cases}
i^2-j+1,&\text{if $i \geq j$},\\
(j-1)^2+i,&\text{if $i < j$.}
\end{cases}$
\\
\\The  sequence $\omega$ is a triangle array read by rows: row number $k$ contains $2k-1$ elements. 
\\
\\Using the angle numeration to the array (4), we have:
\\ ${11},\quad {12},{22}, {21},\quad {13},{23},{33}, {32},{31},\quad {14},{24},{34},{44},{43},{42},{41},\quad ...$
\\
\\The inverse functions are:
\begin{equation} 
i=\text {min}\{t;n-(t-1)^{2}\},\quad \quad
j=\text {min}\{t;t^{2}-n+1\},
\end{equation} 
where $t=\lfloor \sqrt{n-1}\rfloor\ +1$.
\\
\\{\bf Example 5.2.1} \\
Using the angle numeration to the array from Example 2.3.1, we get:
\\$a_{1},\quad a_{k+1},a_{k+2},a_{2},\quad a_{2k+1},a_{2k+2},a_{2k+3},a_{k+3},a_{3},
\\ \quad  a_{3k+1},a_{3k+2},a_{3k+3},a_{3k+4},a_{2k+4},a_{k+4},a_{4}, \quad ...$
\\\\$\omega(n)=a_{m}$, where 
\\$m=(k+1)(\lfloor \sqrt{n-1}\rfloor - \dfrac{|t|+t}{2})+t+ 1 $,\quad $t=n-\lfloor \sqrt{n-1}\rfloor^{2}-\lfloor \sqrt{n-1}\rfloor-1$. 
\\
\\{\bf Example 5.2.2} \\
Using the angle numeration to the array from Example 2.3.2, we have:
\\$a_{1},\quad a_{k+1},a_{k+2},a_{k+1},\quad a_{2k+1},a_{2k+2},a_{2k+3},a_{2k+2},a_{2k+1},
\\ \quad  a_{3k+1},a_{3k+2},a_{3k+3},a_{3k+4},a_{3k+3},a_{3k+2},a_{3k+1}, \quad ...$
\\\\$\omega(n)=a_{m}$, where  $m=(k+1)\lfloor \sqrt{n-1}\rfloor +1 - | n- \lfloor \sqrt{n-1}\rfloor^2 - \lfloor \sqrt{n-1}\rfloor -1 |$.
\\
\\{\bf Example 5.2.3} \\
Using the angle numeration to the array from Example 2.3.3, we obtain:
\\$a_{1},\quad a_{k},a_{1},a_{2},\quad a_{k+1},a_{k},a_{1},a_{2},a_{3},\quad  a_{k+2},a_{k+1},a_{k},a_{1},a_{2},a_{3},a_{4}, \quad ...$
\\
\\The sequence is the generalization of A004739 ($k$=1), A004738 ($k$=2).
\\
\\$\omega(n)=a_{m}$, where $m=kv+(2v-1)(t^{2}-n)+t$,
\\$t=\lfloor \sqrt{n}+1/2 \rfloor+1$, 
$v=\lfloor \dfrac {n-1}{t}\rfloor -t+1$.
\\Using Antonio G. Astudillo formula from A079813, we have formula for $v$. 
\\
\\{\bf Example 5.2.4} \\
Stein proposed the boustrophedonic (``ox-plowing'') option, ``although without giving explicit formula'' [5]. Natural numbers are placed ``by taking clockwise and counter-clockwise turns'' of A081344. 
The formula is below:
\\
\\$n =\begin{cases}
(j-1)^2+j+(-1)^{j-1}(j-i),&\text{if $i \leq j$},\\
(i-1)^2+i-(-1)^{i-1}(i-j),&\text{if $i > j$.}
\end{cases}$
\\
\\The result of the enumeration the array (4) is:
\\ ${11},\quad {12},{22}, {21},\quad {31},{32},{33}, {23},{13},\quad {14},{24},{34},{44},{43},{42},{41},
\\{51},{52},{53},{54}, {55},{45},{35},{25},{15},\quad ...$
\\
\\The property of ``ox-plowing'' is: if $t=\lfloor \sqrt{n-1}\rfloor +1 $ is even, $i, j$ are computed by (5), otherwise  $i, j$ are reversed. Therefore we have formulas for the inverse functions:
\\
\\$i$=Mod$(t;2)$min$\{t;t^{2}-n+1\}$ + Mod$(t+1;2)$min$\{t;n-(t-1)^{2}\}$.
\\$j$=Mod$(t;2)$min$\{t;n-(t-1)^{2}\}$+Mod$(t+1;2)$min$\{t;t^{2}-n+1\}$,
\section{Generalization of Cantor diagonalization}
The Cantor diagonalization assumes numeration by cells 1x1. We consider the cover of the array $\Omega$ by rectangles variable length $l_{m}$ and height $h_{m}$; $l_{m},h_{m} \in \mathbb{Z}^{+} $. The upper right corner is covered by rectangle $l_{1}*h_{1}$, right to it the rectangle $l_{2}*h_{1}$ adjoins, below it the rectangle $l_{1}*h_{2}$ adjoins, and so on. Rectangles with the length $l_{s}$ are placed in the column, rectangles with the height $h_{r}$ are placed in the row. Numeration within all rectangles is row by row, starting from first. The rectangle containing the number $w(i,j)$ is called current. 
\\We define $R$=max$\{r:\sum_{m=1}\limits^rh_{m} < i\}$ and $S$=max$\{s:\sum_{m=1}\limits^sl_{m} < j\}$. 
\\$R,S \in \mathbb{Z}^{*}$. The length and  the height of the current rectangle are $l_{S+1}, h_{R+1}$, respectively.
\\ 
\\Parameter $n$ contains four components: 
\begin{itemize}
\item the number of cells in the stair-step triangle with cathetus $(R+S)$ is 
$\sum_{r=1}\limits^{R+S}h_{r}(\sum_{s=1}\limits^{R+S+1-r}l_{s})$. 
\item the number of cells inside $R$  rectangles, lying above the diagonal passing through the current rectangle is $\sum_{r=1}\limits^{R}h_{r}l_{R+S+2-r}$. 
\item the number of the row inside current rectangle before $w(i,j)$ is $i- \sum_{r=1}\limits^Rh_{r} -1$. 
\item the number of the column $w(i,j)$ inside current rectangle is  $j-\sum_{s=1}\limits^Sl_{s}$. 
\end{itemize}
Therefore
\begin{equation} 
n= \sum_{r=1}\limits^{R+S}h_{r}(\sum_{s=1}\limits^{R+S+1-r}l_{s}) + \sum_{r=1}\limits^{R}h_{r}l_{R+S+2-r} +l_{S+1}(i- \sum_{r=1}\limits^Rh_{r} -1)+ j-\sum_{s=1}\limits^Sl_{s}
\\
\end{equation} 
\\If the  numeration within rectangles is column by column starting from first, formula (6) takes the following form:
\\
\\$n= \sum_{r=1}\limits^{R+S}h_{r}(\sum_{s=1}\limits^{R+S+1-r}l_{s}) + \sum_{r=1}\limits^{R}h_{r}l_{R+S+2-r} +h_{R+1}(j- \sum_{s=1}\limits^Sl_{s} -1) + i-\sum_{r=1}\limits^Rh_{r}$ .
\\
\\Using (2)
\\
\\$t=\dfrac{(R+S)^{2}+3R+S}{2}$ 
\\
\\we get the number of the rectangle. Thus we obtain the formula for numeration within the rectangles row by row or column by column, depending on the parity of $t$:
\\
\\$n= \sum_{r=1}\limits^{R+S}h_{r}(\sum_{s=1}\limits^{R+S+1-r}l_{s}) + \sum_{r=1}\limits^{R}h_{r}l_{R+S+2-r} + 
\\\dfrac{(-1)^{t}+1}{2} \Bigr (l_{S+1}(i- \sum_{r=1}\limits^Rh_{r} -1)+j-\sum_{s=1}\limits^Sl_{s} \Bigl )- 
\\\dfrac{(-1)^{t}-1}{2} \Bigr (h_{R+1}(j- \sum_{s=1}\limits^Sl_{s} -1)+i-\sum_{r=1}\limits^Rh_{r} \Bigl )$.
\\
\\Suppose $t=R+S$; then the above formula is the numeration within the rectangles, depending on the parity of the number of diagonal, formed by the rectangles $l_{s}*h_{r}$.
\\
\\{\bf Example 6.1} \\
We apply the enumeration to the array (4) for rectangles with the lengths of 1,2,3,...  and the heights of 1,2,3,...  Numeration inside rectangles is row by row. The result is:
\\ ${11},\quad {12}, {13}, \quad {21}, {31},  \quad {14},{15},{16}, \quad {22}, {23}, {32},{33},\quad {41},{51},{61},\quad  ...$
\\\\Let the sequences $l_{m}, h_{m}$ be constants: $l_{m} = l, h_{m}=h$ for all $m$. The example $l=2, h=2$ was presented previously without explicit formulas [9]. Using (6), we obviously get:
\\$n= \dfrac{lh}{2}{((R+S)^{2} +3R+S)}+l(i-hR-S-1) +j $.
\\
\\Formulas for $R$ and $S$ take  the following form:
\\$R=\lfloor \dfrac{i-1}{h}\rfloor$ and $S=\lfloor \dfrac{j-1}{l}\rfloor$.
\\
\\{\bf Example 6.2} \\
We apply the enumeration for rectangles with the length of 3 and the height of 2, and the numeration row by row inside it to the array (4). The result is:
\\ ${11}, {12}, {13}, {21}, {22}, {23} \quad {14},{15},{16}, {24}, {25}, {26}, \quad {31},{32},{33},{41},{42},{43}\quad  ...$
\section{Transformation of Sequences as Superposition of Pairing Functions}
Let $f,g$ be pairing functions $\mathbb{Z}^{+}\times\mathbb{Z}^{+}$ $\rightarrow$ $\mathbb{Z}^{+}$, $f \neq g$.  We denote the inverse functions to $f$ by $(f_{1}^{-1}, f_{2}^{-1}) $. We create the array $w(i,j)=a_{f_{1}^{-1}(n), f_{2}^{-1}(n)}$. Using function $g$, we convert the array to the sequence $\omega$ :
$\omega(n)=a_{g(f_{1}^{-1}(n), f_{2}^{-1}(n))}$. 
\\
\\
\\{\bf Example 7.1} \\
Using the enumeration from Example 6.2, we place the sequence $\alpha$, then we form the new sequence by enumeration from Example 5.1.2. The  start of the sequence $\omega$ is:
\\$a_{1},\quad a_{2},a_{4}, \quad a_{5},a_{3},a_{13}, \quad a_{6},a_{14},a_{7}, a_{16},\quad a_{15},a_{10},a_{17},a_{8},a_{32}, \quad 
\\ a_{26},a_{18}, a_{10},a_{33},a_{9}, a_{35},...$,

E-mail:putievskiy@gmail.com

\begin{thebibliography}{9}
\bibitem {Sloane} N. J. A. Sloane, Transformation of integer sequences. \emph{http://oeis.org/transforms.txt/}
\bibitem {Bernstein} M. Bernstein and N. J. A. Sloane, Some Canonical Sequences of Integers, arXiv:0205301v1 [math.CO]
\bibitem{Khovanova} T. Khovanova, How to Create a New Integer Sequence, arXiv:0712.2244 [math.CO] 
\bibitem {Hopcroft}
J.Hopcroft and J.Ullman,  – Introduction to Automata Theory, Languages,
and Computation. – Addison-Wesley, 1979 
\bibitem {Stein} S.Pigeon, Pairing Function, From MathWorld-A Wolfram Web Resource, created by Eric W. Weisstein. http://mathworld.wolfram.com/PairingFunction.html
\bibitem{OEIS} Online Encyclopedia of Integer Sequences (OEIS). \emph{http://oeis.org/}
\bibitem{Kimberling} C. Kimberling, "Numeration systems and fractal sequences," Acta Arithmetica 73 (1995) 103-117.
\bibitem{n-Tuple} D. E. Knuth, The Art of Computer Programming, Volume 4, Fascicle 4, Generating All Trees; History of Combinatorial Generation. - Addison-Wesley, 2006
\bibitem{Pigeon} S. Pigeon, Contributions a la compression de donnees. Ph.d. thesis, Universite de
Montreal, Montreal, 2001
\\
\end{thebibliography}
\end{document}